\documentclass[mcma]{degruyter-journal-a}     

\title{Comparison of Sobol' sequences in financial applications}
\headlinetitle{Comparison of Sobol' sequences}

\lastnameone{Harase}
\firstnameone{Shin}
\nameshortone{S.~Harase}
\addressone{College of Science and Engineering, Ritsumeikan University, 1-1-1 Nojihigashi, Kusatsu, Shiga, 525-8577}
\countryone{Japan}
\emailone{harase@fc.ritsumei.ac.jp}

\lastnametwo{}
\firstnametwo{}
\nameshorttwo{}
\addresstwo{}
\countrytwo{}
\emailtwo{}

\lastnamethree{}
\firstnamethree{}
\nameshortthree{}
\addressthree{}
\countrythree{}
\emailthree{}

\lastnamefour{}
\firstnamefour{}
\nameshortfour{}
\addressfour{}
\countryfour{}
\emailfour{}

\lastnamefive{}
\firstnamefive{}
\nameshortfive{}
\addressfive{}
\countryfive{}
\emailfive{}

\abstract{Sobol' sequences are widely used 
for quasi-Monte Carlo methods that arise in financial applications. 
Sobol' sequences have parameter values called direction numbers, 
which are freely chosen by the user, so there are several implementations of Sobol' sequence generators. 
The aim of this paper is to provide a comparative 
study of (non-commercial) high-dimensional Sobol' sequences by calculating financial models. 
Additionally, we implement the Niederreiter sequence (in base 2) with a slight modification, that is, 
we reorder the rows of the generating matrices, and analyze and compare it with the Sobol' sequences.}



\researchsupported{This work was partially supported by 
 JSPS KAKENHI Grant Numbers JP18K18016, JP26730015, JP26310211, JP15K13460. 
This work was also supported by JST CREST.}

\acknowledgments{}

\dedicated{}

\begin{document}

\section{Introduction}\label{sec:intro}

Monte Carlo (MC) methods are an important numerical tool for pricing many financial derivatives and calculating the Greeks. 
Generally speaking, 
these values can be expressed as mathematical expectations, and the expectations reduce to integrals 
over the $s$-dimensional unit cube $(0, 1)^s$ after a suitable change of variables, 
that is, $\int_{(0,1)^s} f(\mathbf{x}) d \mathbf{x}$ 
for a function $f: (0, 1)^s \to \mathbb{R}$ and $\mathbf{x} := (x_1, \ldots, x_s)$. 
However, it is often difficult to evaluate the exact value analytically 
and the dimension $s$ is over hundreds or thousands, 
so we use MC integration:
\begin{eqnarray} \label{eqn:MC}
\int_{(0,1)^s} f(\mathbf{x}) d \mathbf{x} \approx \frac{1}{N}\sum_{n = 0}^{N-1} f(\mathbf{x}_n),
\end{eqnarray}
where $\{ \mathbf{x}_0, \ldots, \mathbf{x}_{N-1} \} \subset (0, 1)^s$ 
is a point set of independent random samples from 
the uniform distribution on $(0, 1)^s$. 
MC has a probabilistic error of $O(N^{-1/2})$, 
which does not depend on the dimension $s$ but is significantly slow. 
To improve the rate of convergence, 
we apply quasi-Monte Carlo (QMC) methods using {\it low-discrepancy point sets} or {\it sequences}  
that are more uniformly distributed than random points (see \cite{MR2683394,MR1172997} for the precise definition). 
Around the middle of the 1990s, 
a series of studies reported that  
QMC attains a higher rate of convergence than MC for certain types of high-dimensional numerical integration in finance \cite{vACW97a,caflish1997,doi:10.1287/mnsc.42.6.926,doi:10.1080/13504869600000001,vPAS95a}.
Because of this, Sobol' sequences have been widely used since then. 

Sobol' sequences are a class of low-discrepancy sequences originally proposed by Sobol' \cite{MR0219238} in 1967 
and have parameters called {\it direction numbers}, which are freely chosen by the user. 
Thus, there are several implementations of Sobol' sequences with distinct parameter values \cite{Bratley:1988:AIS:42288.214372,MR2001453,MR2429482,Lemieux2004,doi:10.1002/wilm.10056}. 
Some of them have been optimized with the aim of applying them to finance. 
A comparison of Sobol' sequences for high-dimensional problems in finance was presented in \cite{doi:10.1002/wilm.10056}, 
but we want to know further numerical examples, including randomization and effective dimension reduction techniques. 
According to \cite{MR2429482}, Joe and Kuo conducted some preliminary calculations for financial models 
and found that their new Sobol' sequence \cite{MR2429482} 
provided better results in some cases and, at worst, was comparable with the old sequence \cite{MR2001453}; 
however, specific numerical examples were not included in their paper. 

The aim of this paper is to provide a comprehensive comparative 
study of (non-commercial) high-dimensional Sobol' sequences \cite{MR2001453,MR2429482,Lemieux2004} 
in financial applications. 
Niederreiter \cite{MR960233,MR1172997} proposed another class 
of low-discrepancy sequences, called {\it Niederreiter sequences}. 
Recently, Faure and Lemieux \cite{MR3494134} described the relationships between 
Sobol' and Niederreiter sequences in detail. 
Additionally, Faure and Lemieux \cite{FL2018} reported that the Niederreiter sequence (in base $2$)
with a slight modification, i.e., reordering the rows of the generating matrices, demonstrated high performance in some applications. 
Motivated by their report, we also analyze the modified Niederreiter sequence and compare it with Sobol' sequences. 

In the theory of ``analysis of variance'' (ANOVA) decomposition \cite{caflish1997,MR2719643,MR2743492,MR1966664,MR2201179}, 
it is known that the integrand $f(\mathbf{x})$ for certain high-dimensional problems in finance is 
dominated by the first few variables ({\it low effective dimension in the truncation sense})
or 
is well approximated by 
a sum of functions of at most one or two variables ({\it low effective dimension in the superposition sense}), i.e., 
\begin{eqnarray} \label{eqn:ANOVA}
f(\mathbf{x}) = \underset{\mbox{constant}}{\underbrace{f_0}} + \underset{\mbox{order-$1$ terms}}{\underbrace{\sum_{i = 1}^s f_i(x_i)}}
 + \underset{\mbox{order-$2$ terms}}{\underbrace{\sum_{1 \leq i <j \leq s} f_{i, j}(x_i, x_j)}} + (\mbox{small higher-order terms}).
\end{eqnarray}
It is believed that these are reasons why QMC succeeds in high-dimensional numerical integration even 
if the nominal dimension $s$ is over hundreds or thousands. 
The Sobol' sequence provided by Joe and Kuo \cite{MR2429482} was optimized 
so as to have good two-dimensional (2D) projections for the assumption \eqref{eqn:ANOVA}. 
As we shall see later, if the latter condition 
(\ref{eqn:ANOVA}) is satisfied 
but the former condition is not satisfied, that is, 
if $f(\mathbf{x})$ has low effective dimension in the superposition sense 
but high effective dimension in the truncation sense, 
then such optimization seems to be effective. 

The remainder of this paper is organized as follows:  
In Section~\ref{sec:preliminaries}, we review digital nets and sequences, 
the $t$-value, which is a criterion of uniformity, and Sobol' and Niederreiter sequences. 
In Sections~\ref{sec:analysis} and \ref{sec:comparison}, we present our main results. 
In Section~\ref{sec:analysis}, we calculate the frequency of $t$-values 
of Sobol' and Niederreiter sequences for 2D projections in high dimensions
and show that the new Sobol' sequence provided by Joe and Kuo \cite{MR2429482} 
and the modified Niederreiter sequence avoid the existence of extremely large $t$-values. 
In Section~\ref{sec:comparison}, we compare Sobol' and Niederreiter sequences for numerical integration problems, e.g., 
Asian, digital, and basket options, with or without effective dimension reduction. 
In Section~\ref{sec:conclusion}, we conclude this paper. 

\section{Preliminaries} \label{sec:preliminaries}

\subsection{Digital nets and digital sequences} \label{subsec:digital_nets}

Following \cite{MR3038697,MR2683394,MR1172997}, 
we recall a digital method to construct QMC point sets $P$ and (infinite) sequences $S$. 
Sobol' sequences are included in these classes. 
Let  $\mathbb{F}_2 := \{ 0, 1\}$ be the two-element field. 
We perform addition and multiplication over $\mathbb{F}_2$ (or modulo $2$).
\begin{definition}[Digital nets]
Let $s \geq 1$ and $m \geq 1$ be integers.
Let $C_1, . . . , C_s \in \mathbb{F}_2^{m \times m}$ be $m \times m$ matrices over $\mathbb{F}_2$. 
For each $n = 0, \ldots, 2^m-1$, let $n = \sum_{l = 0}^{m-1} n_l 2^{l}$ with $n_l \in \mathbb{F}_2$ be 
the expansion in base $2$. For each $1 \leq i \leq s$, set $(x_{n, i, 0}, \ldots, x_{n, i, m-1})^{\top} := C_i (n_0, \ldots, n_{m-1})^{\top}$, where $\top$ is the transpose, and $x_{n, i} := \sum_{l = 0}^{m-1} x_{n, i, l}2^{-l-1}$. 
Then, the point set $P = \{ \mathbf{x}_n := (x_{n, 1}, \ldots, x_{n, s}) \ | \ n = 0, \ldots, 2^m-1\}$ is called a {\it digital net} over $\mathbb{F}_2$ and $C_1, \ldots,  C_s$ are called the {\it generating matrices} of the digital net $P$.
\end{definition}

The concept of digital nets can be extended to (infinite) sequences $S = \{ \mathbf{x}_0, \mathbf{x}_1, \ldots \} \subset [0, 1)^s$ for $\infty \times \infty$ generating matrices $C_1, \ldots, C_s \in \mathbb{F}_2^{\infty \times \infty}$ and infinite expansions $n = \sum_{l = 0}^{\infty} n_l 2^{l}$ and $x_{n, i} := \sum_{l = 0}^{\infty} x_{n, i, l}2^{-l-1}$ that contain only a finite number of nonzero terms. The resulting sequence $S$ is called a {\it digital sequence}  
and the matrices $C_1, \ldots, C_s$ are called the {\it generating matrices} of the digital sequence. 

\subsection{$(t, m, s)$-nets} \label{subsec:(t, m, s)-nets}

As a quality parameter of uniformity for a point set $P$, 
we recall the definition of the $t$-value. 
See \cite{MR3038697,MR2683394,MR1172997} for details.

\begin{definition}{($(t, m, s)$-nets).}
Let $s \geq 1$, and $t$ be an integer with $0 \leq t \leq m$. 
A point set $P = \{ \mathbf{x}_0, \mathbf{x}_1, \ldots, \mathbf{x}_{2^m - 1}\}$ consisting of $2^m$ points in $[0, 1)^s$ is called a {\it $(t, m, s)$-net} (in base $2$) 
if every subinterval $J = \prod_{i = 1}^s [ {a_i}/{2^{d_i}}, (a_i +1)/{ 2^{d_i}}) \in [0, 1)^s$ 
with integers $d_i \geq 0$ and $0 \leq a_i < 2^{d_i}$ for $1 \leq i \leq s$ and of volume $2^{t-m}$ contains exactly $2^t$ points of $P$. 
\end{definition}
\begin{definition}[$t$-value for a $(t, m, s)$-net]
The minimum $t$ that satisfies the above property is called the {\it $t$-value} for a $(t, m, s)$-net. 
\end{definition}
A point set $P$ is well distributed if the $t$-value is small. 
The integration error is bounded by $O(2^t (\log N)^{s-1}/N)$ for $N= 2^m$ when $f$ is smooth. 
The factor $(\log N)^{s-1}$ is not negligible if $s$ is large, 
but QMC works well for high-dimensional numerical integration in finance 
possibly because $f(\mathbf{x})$ has low effective dimension. 
In the case of digital nets, 
the $t$-value can be easily calculated by some algorithms \cite{MR3085113,Pirsic2001827}. 

\subsection{Sobol' and Niederreiter sequences} \label{subsec:Sobol'}

Sobol' \cite{MR0219238} proposed a construction method for generating matrices 
$C_1, \ldots, C_s \in \mathbb{F}_2^{\infty \times \infty}$ that have a good structure of $(t, m, s)$-nets. 
His sequences are now called {\it Sobol' sequences} and are included in a subclass of {\it generalized Niederreiter sequences} \cite{Tezuka:1993:PAA:169702.169694,rTEZ95a}. 
Recently, Faure and Lemieux \cite{MR3494134} described the relationships between them in detail. 
From this viewpoint, Sobol' sequences are formulated  as follows:
\begin{enumerate}
\item Let $p_1(x) = x \in \mathbb{F}_2[x]$ and $p_i(x) \in \mathbb{F}_2[x]$, $2 \leq i \leq s$, be the $(i-1)$th {\it primitive} polynomials in a list of primitive polynomials that are sorted in non-decreasing order of degree, i.e., $p_2(x) = x+1, p_3(x) = x^2+x+1$, and so on. Let $e_i := \deg ( p_i )$. 
\item For each $1 \leq i \leq s$, set polynomials $g_{i, 0}(x), \ldots, g_{i, {e_i -1}}(x) \in \mathbb{F}_2[x]$ such that 
\begin{eqnarray} \label{eqn:direction number}
\deg g_{i, k} (x)= e_i - 1 -k
\end{eqnarray}
for $0 \leq k \leq e_i-1$, in advance.  These polynomials are the parameters that can be freely chosen by the user 
and correspond one-to one to the so-called {\it direction numbers} (see Remark~\ref{remark:Sobol'} for details). 
\item For $u = 1, 2, \ldots$, consider the formal power series expansion 
\begin{eqnarray} \label{eqn:fps} 
\frac{g_{i, k}(x)}{p_i(x)^u} = \sum_{v = 1}^{\infty} a^{(i)}(u, k, v) x^{-v} \in \mathbb{F}_2((x^{-1})).
\end{eqnarray}
\item Define $C_i = (c_{j, v}^{(i)})_{j \geq 1, v \geq 1} \in \mathbb{F}_2^{\infty \times \infty}$ as $c_{j, v}^{(i)} = a^{(i)}(Q+1, k, v) \in \mathbb{F}_2$ for $1 \leq i \leq s, j \geq 1 , v \geq 1$, 
where 
\begin{eqnarray} \label{eqn:triangular}
j-1 = Q e_i + k,
\end{eqnarray}
with integers $Q = Q(i, j)$ and $k = k(i, j)$ satisfying $0 \leq k \leq e_i-1$. 
Note that each row of $C_i$ corresponds to each formal power series expansion in (\ref{eqn:fps}). 
Note that the conditions (\ref{eqn:direction number}) and (\ref{eqn:triangular}) correspond to the reordering of rows of $C_i$ 
so as to obtain non-singular upper triangular (NUT) matrices. 
\end{enumerate}
The first $2^m$ points $P$ can be viewed as a digital net generated by 
the upper-left $m \times m$ submatrices of $C_1, \ldots, C_s \in \mathbb{F}_2^{\infty \times \infty}$. 
We can prove that $P$ is a $(t, m, s)$-net with the following properties (see \cite{MR2683394, MR1172997}):
\begin{itemize}
\item Each one-dimensional (1D) projection is a $(0, m, 1)$-net, which means that 1D projections have already been optimized, 
that is, each $t$-value is $0$. 
\item The $t$-value is $\leq \sum_{i = 1}^s (e_i-1)$ for any $m$, which means that the initial dimensions have already been optimized. 
\item For any low-dimensional projection, the $t$-value is $\leq \sum(e_i-1)$, 
where $\sum$ is taken over the corresponding projections.
\end{itemize}

The condition (\ref{eqn:direction number}) can be described as  
\begin{eqnarray} \label{eqn:lower}
g_{i, k} (x) = x^{e_i - 1 -k} + (\mbox{lower terms}) \in \mathbb{F}_2[x]
\end{eqnarray}
for $0 \leq k \leq e_i-1$, and a good selection of lower terms makes us obtain $t$-values smaller than those of the above upper bounds. 
\begin{remark} \label{remark:Sobol'}
Sobol' \cite{MR0219238} originally proposed a column-by-column 
construction for generating matrices $C_i$ 
using recurrences of columns based on primitive polynomials $p_i(x)$ for each $1 \leq i \leq s$. 
In this construction, the upper-left $e_i \times e_i$ submatrices of generating matrices $C_i$ are 
initial values, and were originally called the {\it direction numbers}, 
which exactly correspond one-to-one to the polynomials 
$g_{i, 0}(x), \ldots, g_{i, {e_i -1}}(x)$ with (\ref{eqn:direction number}). 
See \cite{MR3494134} for details. 
\end{remark}

\begin{remark} \label{remark:Niederreiter}
Niederreiter \cite{MR960233,MR1172997} proposed another construction 
method for generating matrices $C_1, \ldots, C_s \in \mathbb{F}_b^{\infty \times \infty}$ 
for low-discrepancy sequences, where $b$ is a prime power and $\mathbb{F}_b$ is a finite field with $b$ elements.
These sequences are called {\it Niederreiter} sequences. 
In the standard implementation of \cite{Bratley:1992:ITL:146382.146385} in base $b = 2$, 
the main differences from Sobol' sequences are 
that $p_i(x) \in \mathbb{F}_2[x]$, $1 \leq i \leq s$, are taken to be {\it irreducible} polynomials (sorted in non-decreasing order of degree) instead of primitive polynomials, 
and $g_{i, 0}(x), \ldots, g_{i, {e_i -1}}(x)$ are taken to be
\begin{eqnarray*} \label{eqn:Niederreiter}
g_{i, k}(x) = x^k \in \mathbb{F}_2[x]
\end{eqnarray*}
for $e_i = \deg(p_i)$ and $0 \leq k \leq e_i -1$, instead of polynomials with the condition (\ref{eqn:direction number}). 
Note that there are no freely chosen parameters. 
Note that the resulting generating matrices $C_1, \ldots, C_s \in \mathbb{F}_2^{\infty \times \infty}$ 
do not have the NUT properties (see Figs.~1 and 2 in \cite{MR3494134}), 
so we suffer from the {\it leading-zero phenomenon}, 
that is, there are too many points close to the origin at the beginning of the sequences. 
Additionally, note that we obtain the NUT generating matrices after reordering the rows of $C_i$. 
According to Theorem~4.3 of \cite{MR3494134}, 
such NUT generating matrices are obtained by the original column-by-column construction, 
which implies that the primitivity of $p_i(x)$ is not necessary.
\end{remark}

\section{Analysis of Sobol' and Niederreiter sequences with NUT generating matrices} \label{sec:analysis}

In this section, we compare high-dimensional Sobol' 
and Niederreiter sequences in terms of the $t$-values. 
In 1976, Sobol' and Levitan \cite{SobolLevitan1976} provided direction numbers in terms of Property A and Property A' \cite{SOBOL1976236}, which are the criteria for the equidistribution property of the 1 and 2 most significant bits, respectively. 
In 1988, Bratley and Fox \cite{Bratley:1988:AIS:42288.214372} provided a FORTRAN implementation of the Sobol' sequence using this set of direction numbers up to dimension $40$. 
In our tests, we investigate the following (non-commercial) high-dimensional Sobol' and Niederreiter sequences 
released after 2000:
\begin{enumerate}
\item[(a)] In 2003, Joe and Kuo \cite{MR2001453} provided a Sobol' sequence generator up to dimension $1111$. 
The direction numbers for $1 \leq s \leq 40$ are the same as those of Bratley and Fox \cite{Bratley:1988:AIS:42288.214372}. 
The direction numbers for $40 < s \leq 1111$ are selected so as to satisfy Property A. 
We refer to this generator as Sobol' (JoeKuo03). 
\item[(b)] Lemieux et al.~\cite{Lemieux2004} provided a Sobol' sequence generator up to dimension $360$. 
The direction numbers for $1 \leq s \leq 40$ are the same as those of Bratley and Fox \cite{Bratley:1988:AIS:42288.214372}. 
The direction numbers for $40 < s \leq 360$ are optimized in terms of 
the resolution criterion for eight successive dimensions, see \cite[Chapter~3.5.2]{MR2723077}. 
We use the 2004 version and refer to this generator as Sobol' (Lemieux). 
\item[(c)] In 2008, Joe and Kuo \cite{MR2429482} indicated that the 2003 version of the Sobol' sequence generator has 
a bad structure (i.e., extremely large $t$-values)
for some 2D projections and searched the new direction numbers up to $21201$. 
Their approach was to choose the direction numbers so that (i) Property A holds for $1 \leq s \leq 1111$; and 
(ii) the t-values of 2D projections of the point sets are minimized by proposing the search criterion $D^{(6)}$. 
Consequently, extremely large $t$-values are avoided. 
We refer to this new generator as Sobol' (JoeKuo08). 
\item[(d)] Recently, Faure and Lemieux \cite{MR3494134} 
discussed the Niederreiter sequence (in base $2$) with NUT generating matrices 
after reordering the rows of the generating matrices. 
This sequence can be viewed as the Sobol' sequence based on irreducible polynomials 
$p_1(x) = x, p_2(x), \ldots, p_s(x) \in \mathbb{F}_2[x]$ 
with non-decreasing order of degree and given by $g_{i, k} (x) = x^{e_i - 1 -k}$ 
for $e_i = \deg(p_i)$ and $0 \leq k \leq e_i -1$, so as to satisfy the condition (\ref{eqn:direction number}). 
Additionally, Faure and Lemieux \cite{FL2018} reported that 
the Niederreiter sequence with NUT generating matrices 
already demonstrated high performance in some high-dimensional numerical integrations 
without optimizing the lower terms of $g_{i, k}(x)$ in (\ref{eqn:lower}). 
To confirm their new findings, we implement this generator and refer to it as Niederreiter (NUT). 
The irreducible polynomials and direction numbers are available at \url{https://github.com/sharase/niederreiter-nut}.

\end{enumerate}
In addition, Kucherenko et al.~\cite{doi:10.1002/wilm.10056} released 
commercial software for  Sobol' sequences up to dimension $65536$ 
with Property A for all dimensions and Property A' for successive dimensions, 
but  we exclude it from our tests because it requires a commercial license. 

To compare the sequences (a)--(d), 
we assume that the integrand $f(\mathbf{x})$ satisfies the condition (\ref{eqn:ANOVA}). 
The $t$-values of 1D projections are all $0$, so we calculate the $t$-values of 2D projections. 
Let $m \geq 1$ and 
\[ t(i, s; m) \]
denote the $t$-value of the digital net that corresponds to the $(i, s)$-projection (i.e., the 2D projection of dimensions $i$ and $s$ with $1 \leq i \leq s-1$) of the first $2^m$ points. 
Tables~\ref{table:dimension360} and \ref{table:dimension1024}
show the frequency of all the values of $t(i, s; m)$ $(1 \leq i \leq s-1)$ for $2 \leq s \leq 360$ and $2 \leq s \leq 1024$, respectively. 
Sobol' (Lemieux) is up to dimension $360$, and hence it is excluded in Table~\ref{table:dimension1024}. 
From the tables, there exist extremely large $t$-values for Sobol' (Lemieux) and (JoeKuo03), 
but Niederreiter (NUT) tends to avoid such large $t$-values, as well as Sobol' (JoeKuo08). 
Conversely, the occurrence of small $t$-values (e.g., 1 or 2) for 
Niederreiter (NUT) is more frequent than for Sobol' sequences. 
This implies that Niederreiter (NUT) has high uniformity for 2D projections in high dimensions  
without optimizing direction numbers. 

\newcolumntype{Y}{>{\raggedleft\arraybackslash}p{15pt}} 
\begin{table}[htb] 
\caption{Frequency of $t(i, s; m)$ for $2 \leq s \leq 360$.} \label{table:dimension360}
{\tiny
  \begin{tabular}{|c|c|YYYYYYYYYY|} \hline
 & &  \multicolumn{10}{l|}{Number of occurrences of the $t$-value} \\
$m$ &  & 0 & 1 & 2 & 3 & 4 & 5 & 6 & 7 & 8 & 9 \\ \hline
$10$ & Niederreiter (NUT) & 163 & 11321 & 23097 & 16270 & 7947 & 3495 & 1472 & 576 & 231 & 48 \\
& Sobol' (JoeKuo08) & 214 & 8201 & 20243 & 18004 & 10275 & 4819 & 1924 & 777 & 163 & \\
& Sobol' (Lemieux) & 204 & 8210 & 20040 & 18000 & 10092 & 4865 & 1991 & 863 & 233 & 122 \\
& Sobol' (JoeKuo03) & 204 & 8208 & 20285 & 17854 & 9886 & 4997 & 1961 & 835 & 250 & 140 \\ \hline
$12$ & Niederreiter (NUT) & 71 & 7679 & 22265 & 17458 & 9527 & 4353 & 1954 & 864 & 332 & 79 \\
& Sobol' (JoeKuo08) & 62 & 4752 & 17648 & 19105 & 12303 & 6334 & 2848 & 1127 & 389 & 52 \\ 
& Sobol' (Lemieux) & 49 & 4757 & 17452 & 18973 & 12091 & 6408 & 2842 & 1270 & 498 & 191 \\
& Sobol' (JoeKuo03) & 56 & 4774 & 17342 & 19027 & 12307 & 6249 & 2814 & 1276 & 467 & 215 \\ \hline
$14$ & Niederreiter (NUT) & 21 & 5119 & 19870 & 19164 & 11111 & 5373 & 2375 & 1020 & 405 & 134 \\
& Sobol' (JoeKuo08) & 14 & 2857 & 14942 & 19442 & 14020 & 7581 & 3516 & 1551 & 557 & 140 \\
& Sobol' (Lemieux) & 15 & 2913 & 14912 & 19155 & 13589 & 7531 & 3581 & 1690 & 705 & 329 \\
& Sobol' (JoeKuo03) & 16 & 2864 & 14967 & 19131 & 13669 & 7414 & 3634 & 1698 & 696 & 329 \\ \hline
$16$ &  Niederreiter (NUT) & 6 & 3044 & 16906 & 20478 & 13001 & 6584 & 2903 & 1115 & 434 & 130 \\
& Sobol' (JoeKuo08) & 5 & 1771 & 12568 & 19566 & 14939 & 8566 & 4252 & 1893 & 750 & 285 \\
& Sobol' (Lemieux) & 4 & 1815 & 12586 & 19696 & 14735 & 8257 & 4077 & 1917 & 846 & 408 \\
& Sobol' (JoeKuo03) & 6 & 1745 & 12503 & 19418 & 15039 & 8354 & 4097 & 1902 & 890 & 390 \\ \hline
$18$ &  Niederreiter (NUT) & 8 & 1804 & 14312 & 20990 & 14831 & 7440 & 3267 & 1343 & 455 & 144 \\
& Sobol' (JoeKuo08) & 3 & 1119 & 10985 & 19412 & 15999 & 9175 & 4595 & 2087 & 830 & 335 \\
& Sobol' (Lemieux) & 1 & 1168 & 10897 & 19421 & 15963 & 9064 & 4538 & 2028 & 854 & 380 \\
&  Sobol' (JoeKuo03) & 1 & 1075 & 10866 & 19422 & 15870 & 9236 & 4540 & 2003 & 909 & 397 \\ \hline
$20$ & Niederreiter (NUT) & 3 & 1183 & 12437 & 20942 & 15387 & 8394 & 3806 & 1594 & 636 & 189 \\
& Sobol' (JoeKuo08) & 1 & 787 & 9470 & 19570 & 16932 & 9897 & 4736 & 1992 & 831 & 320 \\
& Sobol' (Lemieux) & 1 & 762 & 9615 & 19188 & 16888 & 9798 & 4884 & 2052 & 868 & 325 \\
&  Sobol' (JoeKuo03) & 1 & 764 & 9526 & 19188 & 17006 & 9834 & 4790 & 2117 & 817 & 351 \\ \hline
  \end{tabular}
  \\*
  \begin{tabular}{|c|c|YYYYYYYYYY|} \hline
 & &  \multicolumn{10}{l|}{Number of occurrences of the $t$-value (continued)} \\
$m$ &  & 10 & 11 & 12 & 13 & 14 & 15 & 16 & 17 & 18 & 19  \\ \hline
$10$ & Niederreiter (NUT) & & & & & & & & & & \\
& Sobol' (JoeKuo08) & & & & & & & & & & \\
& Sobol' (Lemieux) & & & & & & & & & & \\
& Sobol' (JoeKuo03) & & & & & & & & & & \\ \hline
$12$ & Niederreiter (NUT) & 38 & & & & & & & & & \\
& Sobol' (JoeKuo08) & & & & & & & & & & \\ 
& Sobol' (Lemieux) & 56 & 33 & & & & & & & & \\
& Sobol' (JoeKuo03) & 56 & 37 & & & & & & & & \\ \hline
$14$ & Niederreiter (NUT) & 28 & & & & & & & & & \\
& Sobol' (JoeKuo08) & & & & & & & & & & \\
& Sobol' (Lemieux) & 124 & 50 & 17 & 9 & & & & & & \\
& Sobol' (JoeKuo03) & 135 & 50 & 8 & 9 & & & & & &  \\ \hline
$16$ &  Niederreiter (NUT) & 19 & & & & & & & & & \\
& Sobol' (JoeKuo08) & 25 & & & & & & & & & \\
& Sobol' (Lemieux) & 150 & 76 & 28 & 19 & 3 & 3 & & & & \\
& Sobol' (JoeKuo03) & 161 & 71 & 32 & 11 & & 1 & & & & \\ \hline
$18$ &  Niederreiter (NUT) & 24 & 2 & & & & & & & & \\
& Sobol' (JoeKuo08) & 80 & & & & & & & & & \\
& Sobol' (Lemieux) & 167 & 82 & 31 & 17 & 4 & 3 & & 2 & & \\
&  Sobol' (JoeKuo03) & 174 & 73 & 32 & 14 & 7 & 1 & & & & \\ \hline
$20$ & Niederreiter (NUT) & 41 & 8 & & & & & & & & \\
& Sobol' (JoeKuo08) & 82 & 2 & & & & & & & & \\
& Sobol' (Lemieux) & 140 & 55 & 24 & 12 & 4 & 2 & 1 & 1 & & \\
&  Sobol' (JoeKuo03) & 126 & 53 & 32 & 10 & 3 & 2 & & & & \\ \hline
  \end{tabular}
}
\end{table}

\begin{table}[htb] 
\caption{Frequency of $t(i, s; m)$ for $2 \leq s \leq 1024$.} \label{table:dimension1024}
{\tiny
\begin{tabular}{|c|c|YYYYYYYYYY|} \hline
 & &  \multicolumn{10}{l|}{Number of occurrences of the $t$-value} \\
$m$ &  & 0 & 1 & 2 & 3 & 4 & 5 & 6 & 7 & 8 & 9 \\ \hline
$10$ & Niederreiter (NUT) & 1217 & 91368 & 187247 & 131306 & 64096 & 28622 & 12135 & 5079 & 1881 & 825 \\
& Sobol' (JoeKuo08) & 1713 & 66135 & 163425 & 146133 & 81378 & 39763 & 15828 & 6864 & 1920 & 617 \\
& Sobol' (JoeKuo03) & 1761 & 65788 & 163011 & 146129 & 80810 & 40093 & 16015 & 7134 & 2025 & 1010 \\ \hline
$12$ & Niederreiter (NUT) & 358 & 61934 & 178807 & 141504 & 77363 & 35943 & 16214 & 7030 & 2920 & 1188 \\
& Sobol' (JoeKuo08) & 464 & 37931 & 140291 & 154369 & 99570 & 51840 & 23418 & 10403 & 3884 & 1380 \\
& Sobol' (JoeKuo03) & 420 & 37691 & 139255 & 154868 & 99220 & 51931 & 23360 & 10514 & 4002 & 1756 \\ \hline
$14$ & Niederreiter (NUT) & 131 & 41284 & 161893 & 153405 & 87958 & 43456 & 20196 & 9017 & 3904 & 1672 \\
& Sobol' (JoeKuo08) & 113 & 21774 & 116341 & 156401 & 113622 & 62704 & 30350 & 13787 & 5688 & 2326 \\
& Sobol' (JoeKuo03) & 109 & 21745 & 115939 & 155545 & 112483 & 63286 & 30287 & 14127 & 5909 & 2706 \\ \hline
$16$ & Niederreiter (NUT) & 51 & 25934 & 142119 & 159612 & 101194 & 51623 & 24259 & 10993 & 4912 & 2035 \\
& Sobol' (JoeKuo08) & 29 & 12562 & 94650 & 153576 & 123743 & 73026 & 36970 & 17348 & 7471 & 3128 \\
& Sobol' (JoeKuo03) & 37 & 12580 & 94525 & 152964 & 123493 & 72503 & 36228 & 17522 & 7846 & 3583 \\ \hline
$18$ & Niederreiter (NUT) & 20 & 15374 & 118340 & 164947 & 114523 & 60295 & 28790 & 12922 & 5415 & 2118 \\
& Sobol' (JoeKuo08) & 14 & 7362 & 77577 & 148682 & 130582 & 81034 & 42592 & 20643 & 9316 & 3984 \\
& Sobol' (JoeKuo03) & 8 & 7464 & 77049 & 148627 & 131116 & 80131 & 42178 & 20536 & 9189 & 4254 \\ \hline
$20$ & Niederreiter (NUT) & 9 & 8790 & 95347 & 163784 & 127680 & 71672 & 33322 & 14166 & 5747 & 2254 \\
& Sobol' (JoeKuo08) & 4 & 4609 & 64037 & 144019 & 137413 & 87005 & 46421 & 22518 & 10407 & 4664 \\
& Sobol' (JoeKuo03) & 3 & 4524 & 64039 & 143480 & 137250 & 87005 & 46332 & 22273 & 10359 & 4831 \\ \hline
  \end{tabular}
\\*\begin{tabular}{|c|c|YYYYYYYYYY|} \hline

 & &  \multicolumn{10}{l|}{Number of occurrences of the $t$-value (continued)} \\
$m$ &  & 10 & 11 & 12 & 13 & 14 & 15 & 16 & 17 & 18 & 19 \\ \hline
$10$ & Niederreiter (NUT) & & & & & & & & & & \\
& Sobol' (JoeKuo08) & & & & & & & & & & \\
& Sobol' (JoeKuo03) & & & & & & & & & & \\ \hline
$12$ & Niederreiter (NUT) & 432 & 83 & & & & & & & & \\
& Sobol' (JoeKuo08) & 226 & & & & & & & & & \\
& Sobol' (JoeKuo03) & 519 & 240 & & & & & & & & \\ \hline
$14$ & Niederreiter (NUT) & 675 & 139 & 46 & & & & & & & \\
& Sobol' (JoeKuo08) & 622 & 48 & & & & & & & & \\
& Sobol' (JoeKuo03) & 1020 & 436 & 126 & 58 & & & & & & \\ \hline
$16$ & Niederreiter (NUT) & 793 & 220 & 31 & & & & & & & \\
& Sobol' (JoeKuo08) & 1044 & 226 & 3 & & & & & & & \\
& Sobol' (JoeKuo03) & 1439 & 659 & 240 & 118 & 26 & 13 & & & & \\ \hline
$18$ & Niederreiter (NUT) & 767 & 234 & 31 & & & & & & & \\
& Sobol' (JoeKuo08) & 1497 & 446 & 47 & & & & & & & \\
& Sobol' (JoeKuo03) & 1834 & 805 & 331 & 153 & 61 & 32 & 4 & 4 & & \\ \hline
$20$ & Niederreiter (NUT) & 751 & 218 & 33 & 3 & & & & & & \\
& Sobol' (JoeKuo08) & 1940 & 644 & 95 & & & & & & & \\
& Sobol' (JoeKuo03) & 2012 & 942 & 425 & 184 & 70 & 18 & 17 & 9 & 3 & \\ \hline
  \end{tabular}
}
\end{table}

\section{Comparison in financial applications} \label{sec:comparison}

We compare Niederreiter (NUT) and Sobol' sequences in Section~\ref{sec:analysis} from the viewpoint of financial applications. 
In the QMC setting, we apply {\it randomizations} using {\it linear scrambling} 
and {\it digital shift} to point sets (see \cite{MR2000878,MR2519835,MR2723077} for details). 
This technique preserves the $t$-values of $(t, m, s)$-nets and avoids the problem that the first point is always the origin. 
We apply the randomizations $M$ times, make $M$ point sets $\tilde{P}_l = \{ \tilde{\mathbf{x}}_n^{(l)} \} \subset (0,1)^s$ $(l = 1, \ldots, M)$,
and compute $M$ independent estimates of (\ref{eqn:MC}):
\[ Q_l := \frac{1}{N} \sum_{n = 0}^{N-1} f(\tilde{\mathbf{x}}_n^{(l)}). \]
Further, we compute the mean and the standard error of $Q_1, \ldots, Q_M$, i.e., 
\[ \bar{Q} := \frac{1}{M} \sum_{l = 1}^{M} Q_l, \mbox{stderr}(\bar{Q}) := \sqrt{\frac{1}{M(M-1)} \sum_{l = 1}^M (Q_l - \bar{Q})^2}. \]
Throughout this paper, we set $M = 100$ as the number of randomizations. 

\subsection{Asian option}
Assume that under the risk-neutral measure 
the asset price $S_t$ follows the Black--Scholes model (i.e., geometric Brownian motion):
\begin{eqnarray} \label{eqn:black-scholes} 
d S_t = r S_t dt + \sigma S_t d B_t,
\end{eqnarray}
where $r$ is the risk-free interest rate, $\sigma$ is the volatility, $B_t$ is a standard Brownian motion. 
The problem of pricing an Asian call option on the discrete arithmetic average is formulated as follows:   
the payoff function is given by $\max (0, \frac{1}{s} \sum_{i = 1}^s S_{t_i} -K)$, 
where $K$ is the strike price at maturity $T$, 
and a time interval $[0, T]$ is discretized at equally spaced times 
$t_i = i \Delta t$ for $i = 1, \ldots, s$,  where $\Delta t = T/s$. 
Then, the value of the option at time 0 is given by
\begin{eqnarray} \label{eqn:asian}
E \left[ e^{-rT} \max (0, \frac{1}{s} \sum_{i = 1}^s S_{t_i} -K) \right].
\end{eqnarray}
The analytical solution to (\ref{eqn:black-scholes}) is given by $S_t = S_0  \exp ((r - \sigma^2 /2) t + \sigma B_t)$, 
so it is sufficient to simulate sample paths of Brownian motion. 
The standard construction of Brownian motion is to generate $B_{t_i}$ 
sequentially in time: given $B_0 = 0$, 
\begin{eqnarray} \label{eqn:brownian}
B_{t_i} = B_{t_i-1} + \sqrt{\Delta t} Z_{i}, \quad i = 1, \ldots, s,
\end{eqnarray}
where $Z_1, \ldots, Z_s \sim N(0,1)$ are i.i.d.~standard normally distributed random variables. 
The standard construction (\ref{eqn:brownian}) can be written as 
\begin{eqnarray} \label{eqn:cholesky}
(B_{t_1}, \ldots, B_{t_s})^{{\top}} = A (Z_1, \ldots, Z_s)^{\top}, 
\quad A = \sqrt{\Delta t}
\begin{pmatrix}
1 & 0 & \cdots & 0 \\
1 & 1 & \cdots & 0 \\
\vdots & \vdots & \ddots & \vdots \\
1 & 1 & \cdots & 1
\end{pmatrix},
\end{eqnarray}
where $A$ is an $s \times s$ lower triangular matrix. 
Thus, the expectation (\ref{eqn:asian}) can be written as 
\begin{eqnarray*} \label{eqn:asian_interal}
\int_{(0,1)^s} e^{-rT} \max \left( 0, \frac{1}{s} \sum_{i = 1}^{s} S_0 \exp \left[ \left( r - \frac{\sigma^2}{2} \right) t_i + \sigma w_i \right] -K \right) d \mathbf{x},
\end{eqnarray*}
where $\Phi : (0,1) \to \mathbb{R}$ denotes the cumulative distribution function of the standard normal distribution and 
$(w_1, \ldots, w_s)^{\top}  := A(\Phi^{-1}(x_1), \ldots, \Phi^{-1}(x_s))^{\top}$ for $\mathbf{x} = (x_1, \ldots, x_s) \in (0, 1)^s$. 
We use the following parameters: $T = 1, r = 0.1, \sigma = 0.2, S_0 = 100, K = 100$, which were used in \cite{MR1966664}. 

\begin{figure}[htbp]
 \begin{center}
  \includegraphics[width=70mm, angle=270]{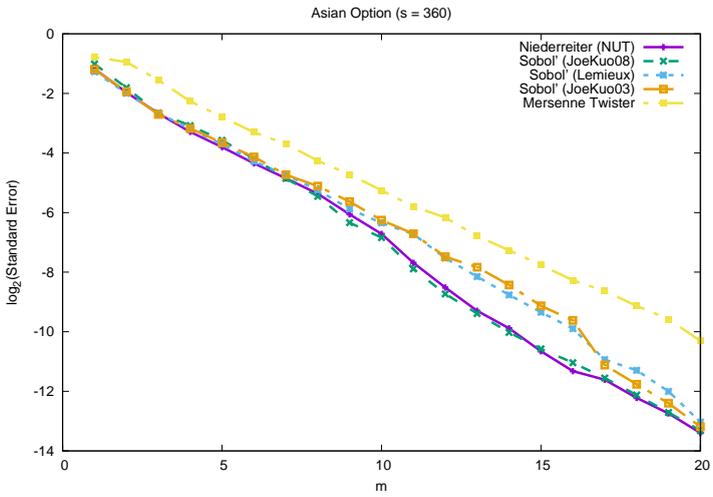}
 \end{center}
 \caption{Comparison of Niederreiter (NUT) and Sobol' sequences for the pricing of an Asian option with $s = 360$.}
 \label{fig:asian360}
\end{figure}

First, we consider the case of dimension $s = 360$. 
Figure~\ref{fig:asian360} shows a summary of the standard error $\mbox{stderr}(\bar{Q})$ in $\log_2$ scale for 
$m = 1, \ldots, 20$. 
In our experiments, 
we applied QMC methods based on Sobol' and Niederreiter sequences (a)--(d) 
and crude MC methods using random number sequences from Mersenne Twister \cite{Matsumoto:1998:MTE:272991.272995}. 
For this, we observed that Niederreiter (NUT) and Sobol' (JoeKuo08) are more effective than the others, particularly for $m = 10, \ldots, 18$, which are often used in practice.  
Additionally, our result seems to agree with the frequency of $t$-values for 2D projections in Table~\ref{table:dimension360}.  
Thus, it is inferred that the Asian option using the standard construction (\ref{eqn:brownian}) 
has low effective dimension in the superposition sense 
but high effective dimension in the truncation sense.

\begin{figure}[htbp]
 \begin{center}
  \includegraphics[width=70mm, angle=270]{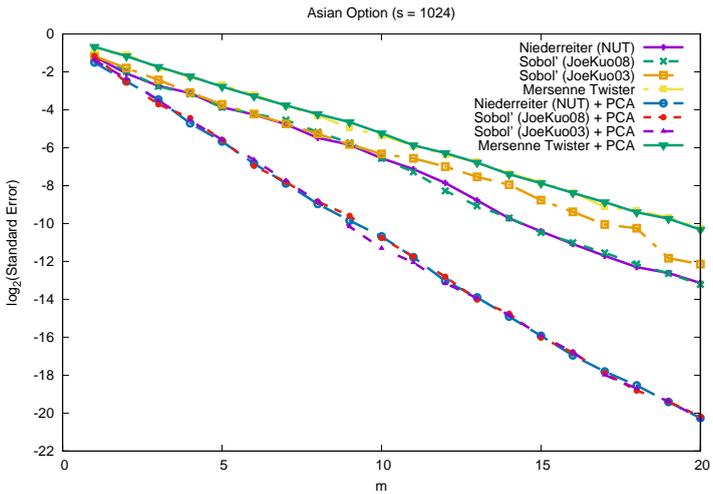}
 \end{center}
 \caption{Comparison of Niederreiter (NUT) and Sobol' sequences for the pricing of an Asian option using the standard and PCA constructions for $s = 1024$.}
 \label{fig:asian1024}
\end{figure}

Further, we consider the higher dimensional case $s = 1024$. 
Figure~\ref{fig:asian1024} shows a summary of $\mbox{stderr}(\bar{Q})$ in $\log_2$ scale. 
In the standard construction (\ref{eqn:brownian}), Niederreiter (NUT) and Sobol' (JoeKuo08) also provide better results than Sobol' (JoeKuo03). 
Further, we recall {\it dimension reduction techniques}, such as the principal component analysis (PCA) 
construction \cite{vACW97a} for generating Brownian motion $B_t$, 
which enhance the efficiency of QMC methods. 
Here, the sampled Brownian motion $(B_{t_1}, \ldots, B_{t_s})^{\top}$ is normally distributed with mean $\bf{0}$ and covariance matrix 
$C = ( \min (t_i, t_j) )_{i, j = 1}^s$, i.e., $(B_{t_1}, \ldots, B_{t_s})^{{\top}} \sim N(\mathbf{0}, C)$. 
Generally, we obtain the equivalent paths of Brownian motion:
\begin{eqnarray*} \label{eqn:decomposition}
(B_{t_1}, \ldots, B_{t_s})^{\top} = A (Z_1, \ldots, Z_s)^{\top}, \quad (Z_1, \ldots, Z_s)^{\top} 
\sim N( \mathbf{0}, I_s), 
\end{eqnarray*}
provided we apply the change of variables $\mathbf{x} = A \mathbf{z}$ with $A A^{\top} = C$.    
The matrix $A$ in the standard construction (\ref{eqn:cholesky}) is the Cholesky matrix of $ C$, i.e., $A A^{\top} = C$. 
Conversely, the PCA construction is a method used to choose $A = [ \sqrt{\lambda_1} \mathbf{v}_1,\cdots, \sqrt{\lambda_s} \mathbf{v}_s]$, where $\lambda_1 \geq \cdots \geq \lambda_s$ are the eigenvalues 
and $\mathbf{v}_1, \ldots, \mathbf{v}_s$ are the corresponding unit-length eigenvectors of $C$. 
In our experiment, the Niederreiter (NUT) and Sobol' sequences with PCA outperform those with the standard construction, 
but those with PCA have exactly the same convergence rates. 
Our result implies that PCA transforms the integrand 
so as to have low effective dimension in the truncation sense, that is, 
the important variables are concentrated in the first few dimensions (e.g., $\leq 2$ or $3$). 
However, for Sobol' sequences,
the lower terms of $g_{i, k}(x)$ in (\ref{eqn:lower})
have almost no choice and are fixed in these first dimensions because the degree $e_i$ in (\ref{eqn:lower}) is sufficiently small. 
Thus, it seems to be difficult to expect further improvement for dimension reduction techniques 
as a result of changing the direction numbers for Sobol' sequences. 

We also tested the Brownian bridge (BB) construction \cite{MR1398000} as another dimension reduction technique and 
observed that there is no difference among Niederreiter (NUT) and Sobol' sequences for the convergence rates, 
which are better than the standard construction but worse than PCA, 
so we omitted the results. 

\subsection{Digital option}

\begin{figure}[htbp]
 \begin{center}
  \includegraphics[width=70mm, angle=270]{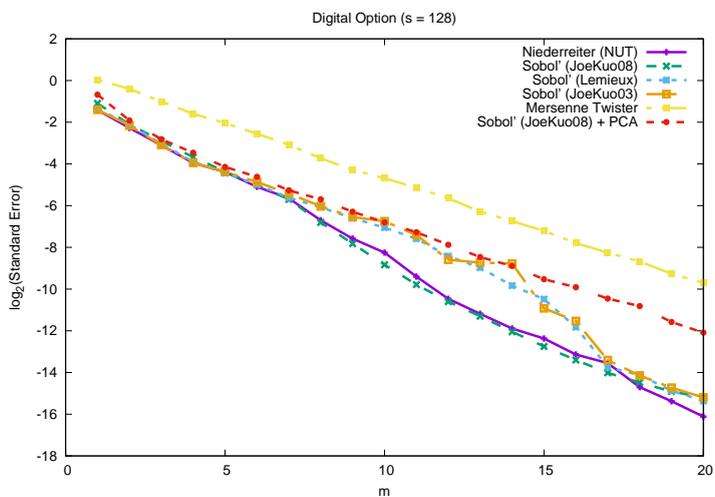}
 \end{center}
 \caption{Comparison of Niederreiter (NUT) and Sobol' sequences for the pricing of a digital option with $s = 128$.}
 \label{fig:digital128}
\end{figure}

Assume that the asset price $S_t$ follows the Black--Scholes model (\ref{eqn:black-scholes}). 
Papageorgiou~\cite{PAPAGEORGIOU2002171} considered the following {\it digital option}:
\begin{eqnarray} \label{eqn:digital}
E \left[ \frac{1}{s} \sum_{i = 1}^s (S_{t_i} - S_{t_{i-1}})_{+}^{0}S_{t_i} \right],
\end{eqnarray}
where $(x)_{+}^{0}$ is equal to 1 if $x > 0$ and is $0$ otherwise, $x \in \mathbb{R}$. 
He indicated that effective dimension reduction techniques perform worse than the standard construction (\ref{eqn:brownian}). 
Wang and Tan \cite{doi:10.1287/mnsc.1120.1568} and Wang \cite{MR3500605} 
found that if the paths are generated by the standard construction, 
then the discontinuities of the payoff function of the sum of the indicator functions 
are aligned with the coordinate axes, so good performance is expected, 
but BB and PCA do not have this type of discontinuity. 
Thus, the standard construction is a good choice in this QMC setting. 
Figure~\ref{fig:digital128} shows a summary of $\mbox{stderr}(\bar{Q})$ in $\log_2$ scale. 
We used the parameters $s = 128, T = 1, r = 0.045, \sigma = 0.3, S_0 = 100$ from \cite{PAPAGEORGIOU2002171}. 
Indeed, PCA is worse than the standard construction. 
Niederreiter (NUT) and Sobol' (JoeKuo08) with good 2D projections are useful for such a problem. 
Note that this example is very simple and the value (\ref{eqn:digital}) can be calculated analytically.

\subsection{Basket option}

\begin{figure}[htbp]
 \begin{center}
  \includegraphics[width=70mm, angle=270]{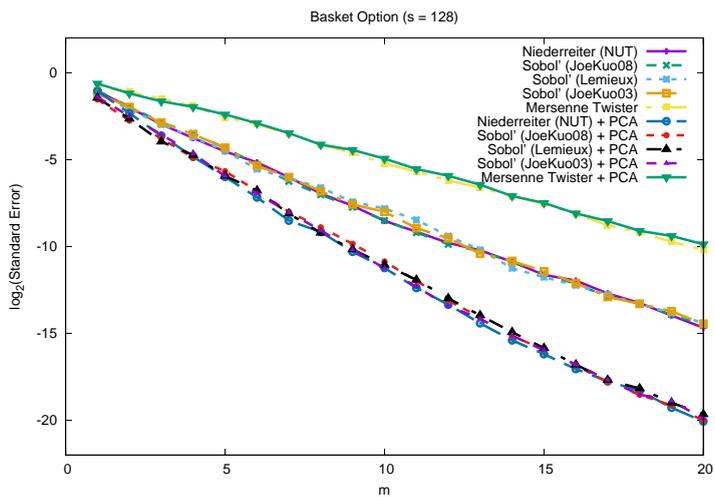}
 \end{center}
 \caption{Comparison of Niederreiter (NUT) and Sobol' sequences for the pricing of a basket option with number of assets $s = 128$.}
 \label{fig:basket128}
\end{figure}

Following \cite{MR1966664,MR2201179}, under the risk-neutral measure, 
we consider a European-style basket call option on the arithmetic average over $s$ 
assets $S_t^{(1)}, \ldots, S_t^{(s)}$, and assume 
that each asset satisfies 
\begin{eqnarray} \label{eqn:basketSDE}
d S_t^{(i)} = r S_t^{(i)} dt + \sigma^{(i)} S_t^{(i)} d B_t^{(i)} \quad (i =1 , \ldots, s),
\end{eqnarray}
for a mean return parameter $r$ and volatility parameters $\sigma^{(i)}$. 
Assume $B_t^{(1)}, \ldots, B_t^{(s)}$ are correlated Brownian motions with correlations $\rho_{ij}$, 
and the terminal pay off at $T$ is given by $\max ( 0, \frac{1}{s} \sum_{i = 1}^s S_T^{(i)} -K )$. 
For this, we compute the price of the basket option:
\begin{eqnarray} \label{eqn:basket}
E \left[ e^{-rT}  \max \left( 0, \frac{1}{s}\sum_{i = 1}^s S_T^{(i)} -K \right) \right].
\end{eqnarray}
Note that $s$ is the number of assets, not the number of discretization steps. 
The solutions to (\ref{eqn:basketSDE}) are given by $S_t^{(i)} = S_0^{(i)} \exp ((r - (\sigma^{(i)})^2 /2) t + \sigma^{(i)} B_t^{(i)})$.
Here, the random vector $(B_T^{(1)}, \ldots, B_T^{(s)})^{\top}$ is normally distributed with mean $\mathbf{0}$ and 
covariance matrix $C = (\rho_{ij} T)_{i, j = 1}^s$. 
Let $(Z_1, \ldots, Z_s)^{\top} \sim N( \mathbf{0}, I_s)$. 
The standard construction for generating Brownian motion is $(B_T^{(1)}, \ldots, B_T^{(s)})^{\top} = A (Z_1, \dots, Z_s)^{\top}$, 
where $A$ is the Cholesky matrix of $C$. 
By contrast, the PCA chooses $A = [ \sqrt{\lambda_1} \mathbf{v}_1, \cdots, \sqrt{\lambda_s} \mathbf{v}_s]$, where $\lambda_1 \geq \cdots \geq \lambda_s$ are the eigenvalues 
and $\mathbf{v}_1, \ldots, \mathbf{v}_s$ are the corresponding unit-length eigenvectors of $C$. 
Expectation (\ref{eqn:basket}) is expressed as 
\begin{eqnarray*} \label{eqn:basket_interal}
\int_{(0,1)^s} e^{-rT} \max \left( 0, \frac{1}{s} \sum_{i = 1}^{s} S_0^{(i)} \exp \left[ \left( r - \frac{(\sigma^{(i)})^2}{2} \right) T + \sigma^{(i)} w_i \right] -K \right) d \mathbf{x},
\end{eqnarray*}
where $(w_1, \ldots, w_s)^{\top}  := A(\Phi^{-1}(x_1), \ldots, \Phi^{-1}(x_s))^{\top}$.
We set the parameters $s = 128, T = 1, r = 0.1, \sigma^{(i)} = 0.2, \rho_{ij} = 0.3 (i \neq j),  S_0^{(i)} = 100, K = 100$, 
which are taken from \cite{MR1966664,MR2201179}, 
and conduct experiments on the standard and PCA constructions. 
Figure~\ref{fig:basket128} shows a summary of $\mbox{stderr}(\bar{Q})$ in $\log_2$ scale. 
Unlike the previous examples, the Niederreiter (NUT) and Sobol' sequences 
using the standard (Cholesky) construction have exactly the same convergence rates. 
According to \cite[Table~3 and 6]{MR1966664} and \cite[Table~3 and Table 5]{MR2201179}, 
it is inferred that the value of basket options using the standard (Cholesky) construction 
is determined by depending on the first few variables or depending on a high proportion to order-1 terms 
$\sum_{i = 1}^s f_i(x_i)$ in (\ref{eqn:ANOVA}), compared with those of Asian options. 

\subsection{Asian option under the Heston model}

As a more complicated model, 
under the risk-neutral measure, 
we consider the pricing of an Asian call option (\ref{eqn:asian}) 
with maturity $T$ and strike $K$ written on an asset whose price
process $S_t$ satisfies the Heston stochastic volatility model:
\begin{eqnarray*}
d S_t  &=& r S_t dt + \sigma_t S_t \left[ \rho d B_t^{(1)} + \sqrt{1-\rho^2} dB_t^{(2)} \right],\\ 
d \sigma_t^2 & = & \kappa \left[ \theta - \sigma_t^2 \right] dt + \xi \sigma_t dB_t^{(1)},
\end{eqnarray*}
where $\sigma_t^2$ is the volatility process, $B_t^{(1)}$ and $B_t^{(2)}$ are two independent standard Brownian motions, 
$r$ is the risk-free interest rate, 
$\kappa$ is the {\it speed of mean reversion}, $\theta >0$ is the {\it long-run mean variance}, 
$\xi$ is the {\it volatility of the volatility}, $\rho$ is the correlation between the Brownian motions driving 
$S_t$ and $\sigma_t^2$. 
The volatility process $\sigma_t^2$ follows a CIR process, which is always positive under the assumption $2 \kappa \theta > \xi^2$.  
We use the Euler--Maruyama scheme with $s$ steps to discretize both $S_t$ and $\sigma_t^2$ as in \cite[Fig.~7.3 in Chapter~7.2.1]{MR2723077}. Let $\Delta t = T/s$. 
Then, we need $2s$-dimensional points to simulate both $S_{t_i}$ and $\sigma_{t_i}^2$ for $t_i = i \Delta t$ $(i = 1, \ldots, s)$. 

\begin{figure}[htbp]
 \begin{center}
  \includegraphics[width=70mm, angle=270]{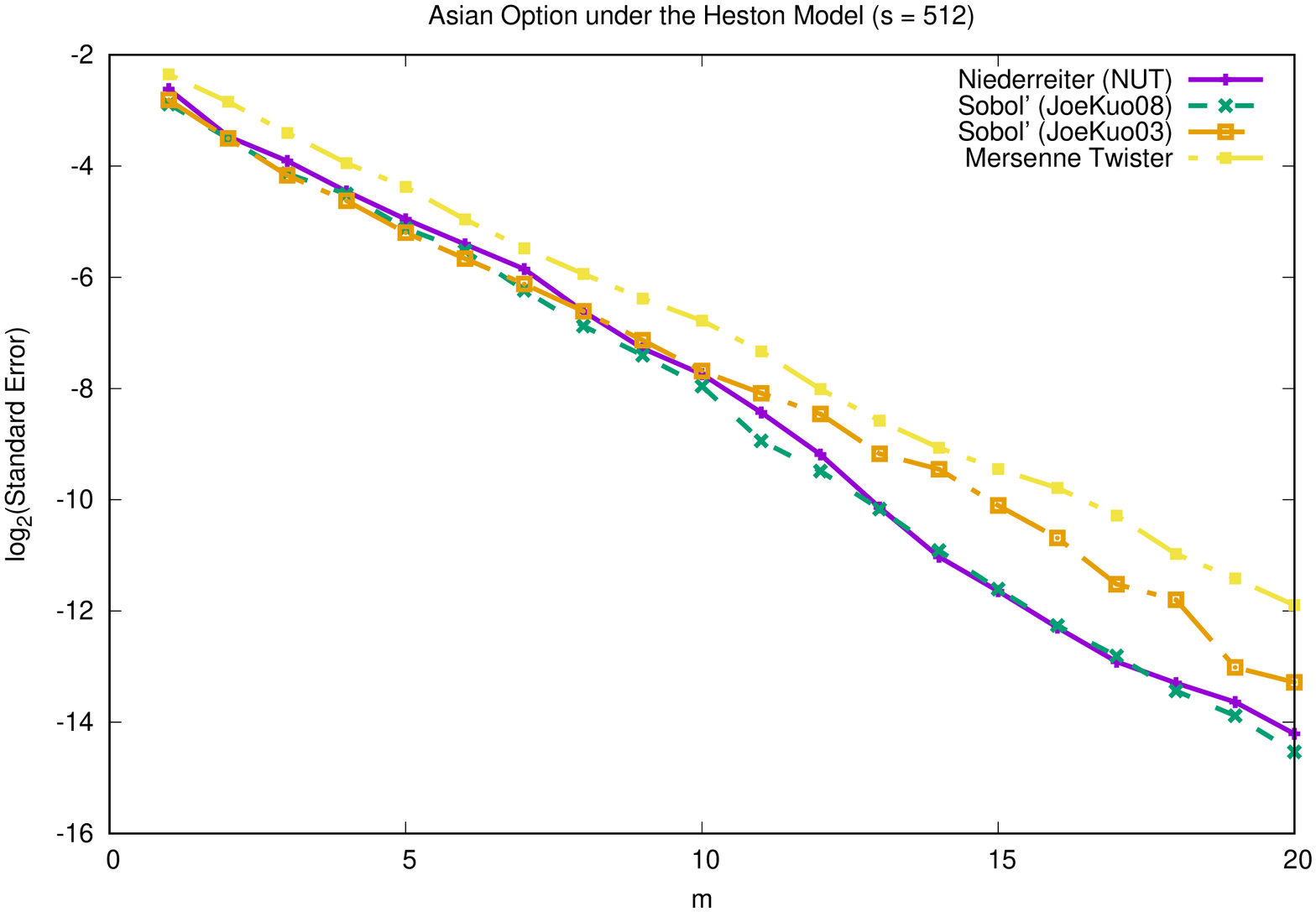}
 \end{center}
 \caption{Comparison of Niederreiter (NUT) and Sobol' sequences for the pricing of an Asian option under the Heston model for $s = 512$.}
 \label{fig:heston512}
\end{figure}

Figure~\ref{fig:heston512} gives results for an Asian option under the Heston model with $s = 512$. 
We use the parameters $T = 0.5, r = 0, \kappa = 2, \theta = 0.01, \xi = 0.1, \rho = 0.5, S_0 = 100, \sigma_0 = 0.1, K = 100$, which are from \cite[Chapter~7.3]{MR2723077}. 
Note that Niederreiter (NUT) and Sobol' (JoeKuo08) give better results than Sobol' (JoeKuo03).  


\section{Concluding remarks} \label{sec:conclusion}

Sobol' sequences have been used successfully in high-dimensional numerical integration in financial applications. 
There are several implementations of Sobol' sequence generators with distinct direction numbers, 
so it is natural to assess which of them is better. 
Hence, we tested Sobol' sequences for calculating financial models and 
observed that the Sobol' sequence with good 2D projections \cite{MR2429482}
outperforms the previous Sobol' sequences \cite{MR2001453,Lemieux2004}, 
particularly in the case in which the integrands have low effective dimension in the superposition sense 
but high effective dimension in the truncation sense. 
Additionally, we implemented the Niederreiter sequence with NUT generating matrices, 
suggested by Faure and Lemieux \cite{MR3494134,FL2018}. 
Surprisingly, the modified Niederreiter sequence has already had good 2D projections, 
and had high performance for calculating some financial models without optimizing direction numbers. 

Finally, we mention the possibility of the further improvement of Sobol' sequences. 
In fact, we attempted the further improvement of Sobol' type digital nets by optimizing direction numbers 
in terms of the framework of generalized Niederreiter sequences based on irreducible polynomials. 
For example, we attempted to search direction numbers so as to have even better $t$-values for 2D projections 
or better three-dimensional projections in addition to good 2D projections. 
The aim appears to be theoretically achieved, but we could not observe a clear difference from 
the Niederreiter sequence with the NUT generating matrices for calculating financial models. 
Perhaps, there might be a limit to the further improvement of Sobol' sequences 
by optimizing direction numbers from a practical perspective. 
Therefore, further improvement of Sobol' sequences is left for future work.

\bibliography{sobolbib}
\end{document}